\documentclass{article}

\usepackage{amsmath,amssymb}
\usepackage{array,booktabs,ragged2e,caption} 
\usepackage{graphicx}
\usepackage[colorinlistoftodos,prependcaption,textsize=tiny]{todonotes}
\usepackage{zibtitlepage}
\usepackage{url}
\usepackage[a4paper,left=.20\paperwidth,right=0.20\paperwidth]{geometry}

\begin{document}
\ZTPAuthor{Felix Hennings}
\ZTPTitle{Benefits and Limitations of Simplified Transient Gas Flow Formulations}
\ZTPNumber{17-39}
\ZTPMonth{Juli}
\ZTPYear{2017}

\zibtitlepage

\title{Benefits and Limitations of Simplified Transient Gas Flow Formulations}
\author{Felix Hennings}

\maketitle

\begin{abstract}
Although intensively studied in recent years, the optimization of the transient (time-dependent) control of large real-world gas networks is still out of reach for current state-of-the-art approaches.
For this reason, we present further simplifications of the commonly used model, which lead to a linear description of the gas flow on pipelines.
In an empirical analysis of real-world data, we investigate the properties of the involved quantities and evaluate the errors made by our simplification.
\end{abstract}

\section{Introduction}
For the past years, the mathematics of gas transport have been intensively studied, mainly focusing on the stationary (time-independent) case, which can be applied to planning scenarios for example~\cite{KocHilPfeSch2015}\cite{PfeFueGeiGei2014}.
However, when aiming to optimize the actual short-term control of the technical gas network elements, we have to consider the time dependent so-called transient case.
Here, research is still in the early stages and current state-of-the-art approaches cannot solve instances of large real-world network size~\cite{RBor2015}.

One difficulty are the Euler Equations~\cite{KocHilPfeSch2015} describing the one dimensional gas flow in a cylindric pipeline, a set of nonlinear hyperbolic partial differential equations. For the isothermal case they can be stated as
\begin{align*}
    \frac{\partial \rho}{\partial t} + \frac{\partial (\rho\,v)}{\partial x} &=0 \\
    \frac{\partial (\rho\,v)}{\partial t} + \frac{\partial (p+\rho\,v^2)}{\partial x} + \lambda\frac{|v|\,v}{2D}\rho + g\,\rho\,h' &=0,
\end{align*}
where $x$ denotes the position in the pipe, $t$ the time, $\rho$ the density of the gas, $v$ the velocity of the gas, $p$ the pressure of the gas, $\lambda$ the friction factor of the pipeline, $D$ the diameter of the pipeline, $g$ the gravitational acceleration and $h'$ the constant slope of the pipe.
Note that $\rho$, $v$ and $p$ depend on $x$ and $t$.
The second equation can be further simplified by assuming the terms $\partial_t (\rho\,v)$ and $\partial_x (\rho\,v^2)$ to be small as in~\cite{EhrSte2003}.
We finally reduce the number of variables by rewriting the constraints using the equation of state for real gases $p = R_s\,\rho\,T\,z(p)$ and the definition of the mass flow $q = A \rho v$ with $A = D^2 \pi / 4$ being the cross sectional area of the pipe as
\begin{align}
    \frac{A}{R_s\,T\,z}\frac{\partial p}{\partial t} + \frac{\partial q}{\partial x} &=0 \label{eq:final_fst} \\
    \frac{\partial p}{\partial x} + \frac{\lambda\,R_s\,T\,z}{2DA^2}\frac{|q|\,q}{p} + g\,h'\,\frac{p}{R_sTz} &=0. \label{eq:final_snd}
\end{align}
Here $R_s$ denotes the specific gas constant, $T$ the constant gas temperature and $z(p,T)$ the compressibility factor, which is often assumed to be constant and hence just stated as $z$.

This model of gas flow in pipelines still contains non-convex terms, which introduce a lot of complexity to any model aiming to solve these equations.
For this reason, we will present an additional simplification of the constraints and investigate the resulting theoretical properties as well as evaluate the caused errors based on historic flow data of real pipelines.

\section{A linearization approach}
\label{sec:linearization}
The non-convexity of the stated equations is based in the second term in \eqref{eq:final_snd}
\begin{equation*}
   f := \frac{\lambda\,R_s\,T\,z}{2DA^2}\frac{|q|\,q}{p}
\end{equation*}
describing the \emph{friction-based pressure difference per meter} on a pipeline.
Using the definition of mass flow and the equation of state above we get a definition of the velocity in terms of pressure and mass flow as
\begin{equation}
   v = \frac{R_s\,\,T\,z}{A}\frac{q}{p} \quad\Rightarrow\quad |v|=\frac{R_s\,\,T\,z}{A}\frac{|q|}{p}. \label{eq:velo}
\end{equation}
We can now rewrite \eqref{eq:final_snd} as
\begin{equation*}
   \frac{\partial p}{\partial x} + \frac{\lambda\,|v|}{2DA}q + g\,h'\,\frac{p}{R_sTz} =0
\end{equation*}
and observe that the equation becomes linear if we assume the absolute velocity in the friction term to be constant, that is $|v|=v_c$.
Note that we do not restrict ourselves to one flow direction, since we fix the \emph{absolute} velocity.
Furthermore, we only fix the absolute velocity in the friction term.
Hence, the \emph{actual} absolute velocity calculated from $p$ and $q$ might be different from $v_c$.

If the proposed simplification can be verified to be reasonable, the overall modeling complexity would decrease drastically.
However, since the friction-based pressure difference scales linearly with the velocity, the assumption of fixed velocity might easily lead to large errors in terms of pressure differences.
On the other hand, both the friction-induced pressure difference and the absolute velocity increase with increasing absolute flow values.
As a consequence, overestimating the absolute velocity should in general be more favorable for minimizing the error in terms of pressure differences.

\section{Analysis of real-world data}
\label{sec:real_data_analysis}
In order to see if the approximation of the friction term by using a constant absolute velocity is reasonably close to the actual friction-induced pressure differences, we use real pipeline data in the network of our project partner OGE~\cite{OGE}, which is the biggest gas network operator in Germany.
For this network, a history of states is given, measured every three minutes over a period of two years.
There are two types of pipelines we consider here: (a)~four large pipes $A$ to $D$, which are used to transport gas between large network intersection areas, and (b)~two small pipelines $E$ and $F$, which are part of the network section connecting customers nodes with bigger pipelines.
An overview of the properties of the six pipes can be found in Table~\ref{table:pipe_data}.

{
\small
\centering
\captionsetup{justification=centering}
\captionof{table}{Properties of the analyzed pipelines.
                  The averages are computed over two years.
                  For the pressure averages a stationary formula (Lem~2.3 from~\cite{KocHilPfeSch2015}) is used.
                  The calculation of the velocities \eqref{eq:velo} uses for the compressibility factor $z$ the formula of Papay~\cite{Pap1968}\cite{Sal2002} and as an aggregated gas mixture computed from the mixtures at entries using a formula for mixtures at junctions from~\cite{KocHilPfeSch2015} Chap~2.
                  The last column denotes the percentage of times when gas was flowing into the main direction.}
\label{table:pipe_data}
\begin{tabular}{lrrrrr}
Pipe & $L$ [km] & $D$ [mm] & avg $p$ [bar] & avg $|v|$ [m/s] & flow in main direction \\
\midrule
A    & 16 & 1000 & 56 & 4.2 & 100 \%\\
B    & 16 &  900 & 63 & 3.4 &  99 \%\\
C    & 15 & 1100 & 70 & 2.5 &  99 \%\\
D    & 20 & 1100 & 71 & 1.4 &  76 \%\\
E    &  3 &  400 & 54 & 2.7 &  93 \%\\
F    &  2 &  300 & 16 & 4.4 & 100 \%
\end{tabular}

\smallskip
}

One idea to choose a fixed absolute velocity value for each pipe is to set this value to a mean velocity computed over a given time period.
The fixed velocity error would be small if a major part of absolute velocity values were to lie within a small interval and hence near the mean value.
To investigate this, we plotted the absolute velocity values occurring on all pipelines over the two years in Figure~\ref{fig:velos}.

We can observe that the pipes $A$, $B$ and $C$ share a distinct characteristic with a high population of similar values and steep tails of the distributions.
For example, for pipe $A$, the absolute difference between the 10th percentile and the 90th percentile is $3.10$.
Thus, only a relative error of less than $36\%$ has to be taken into account for these values, when assuming a mean value of $4.290$ m/s.
For $B$ and $C$ the same percentiles yield a relative error of less than $44\%$ resp. $54\%$ for mean values of $3.445$ resp. $2.610$ m/s.

In contrast, the other pipes have a much larger span of velocity values, even when ignoring the first $13\%$ of values smaller than $0.02$ m/s of pipe $E$.
We can conclude that for the bigger pipes with unique flow directions, namely $A$, $B$ and $C$, a precalculated constant mean velocity value should lead to relatively small errors in terms of friction loss.

A second idea to approximate the velocity term, would be to fix the velocity to some value known from the recent history, such as the velocity value given in the initial state of the gas network control problem for example.
To estimate a possible error here, we compute for each time step the relative velocity changes over one week and average these over the whole time period of two years.
The results are given in Figure~\ref{fig:velochanges}.
In our discussion we focus on the value at 48 hours since this is a typical time horizon for short-term gas network control problems.

{
    \small
    \centering
    \captionsetup{justification=centering}
    \begin{minipage}[t]{0.48\textwidth}
        \centering
        \includegraphics[width=\textwidth]{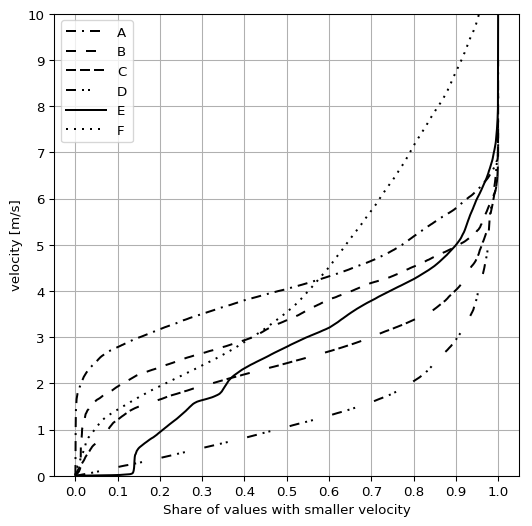}
        \captionof{figure}{Sorted absolute velocity values of each pipe over two years (cumulated distribution functions)}
        \label{fig:velos}
    \end{minipage}
    \hspace*{0.04\textwidth}
    \begin{minipage}[t]{0.48\textwidth}
        \centering
        \includegraphics[width=\textwidth]{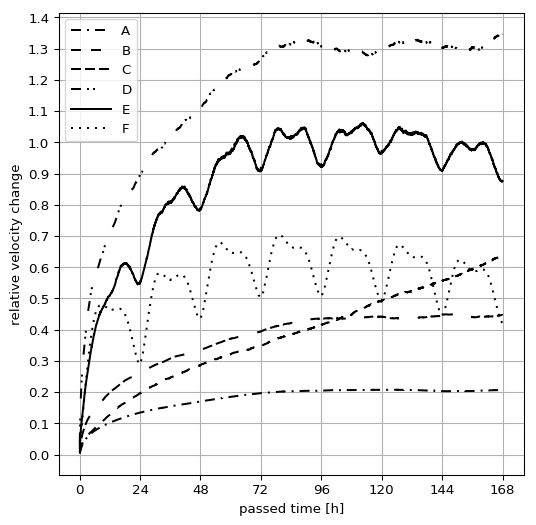}
        \captionof{figure}{Relative velocity change over time. Average values over two years ignoring abs. velocities below $0.02$ m/s.}
        \label{fig:velochanges}
    \end{minipage}
}

In the picture we see, that for the pipes $A$, $B$ and $C$ the velocities change only slowly over short time periods, i.e. on average less than $0.35$ m/s in 48 hours for each of the pipes.
Pipe $A$ and $C$ even seem to reach some constant level of relative velocity change over time.
For the other three pipes we have higher relative velocity changes in the first 48 hours, especially for pipe $D$.
However, the smaller pipes $E$ and $F$ seem to have some daily pattern with a local minimum every 24 hours, which should be taken into account when substituting the velocity with historical values.

\section{Determining fixed velocity values}
After studying the historical velocity values on the six pipes, we will now compute concrete predefined velocity values and examine the actual errors in terms of pressure difference.
We use the two approaches already briefly mentioned above:
Calculating a constant velocity for each pipe based on a large set of historical data (approach $\mathcal{A}$) and taking the velocity, which has been measured on the pipe exactly 48 hours before (approach $\mathcal{B}$).

The constant velocity of $\mathcal{A}$ is calculated as the velocity that minimizes the sum over time of squared pressure difference errors on the whole pipe length $fL$.
These errors are derived based on \eqref{eq:final_snd} by replacing the derivative with the corresponding difference quotient and using $v_c$.
Again we use Papay for the compressibility factor $z$ and the formula of Nikuradse~\cite{Nik1950} for the friction factor $\lambda$.
To have an unbiased evaluation, we use only the first year of the time period for calculating the constant velocity and compare the results of the two approaches on the basis of the second year data.

{
\centering
\small
\captionsetup{justification=centering}
\captionof{table}{Results of approach $\mathcal{A}$ including the calculated constant absolute velocity $v_c$, the average pressure difference error, the corresponding maximum error, the average and maximum real friction-based pressure differences and the ratio of errors to friction values for the average and the maximum case.}
\label{table:results_A}
\begin{tabular}{lrrrrrrr}
Pipe & $v_c$ & avg err & max err & avg $|fL|$ & max $|fL|$ & avg err/ & max err/\\
 & in m/s & in bar & in bar & in bar & in bar & avg $|fL|$ & max $|fL|$\\
\midrule
A &    4.969 &    0.130 &    0.801 &    0.621 &    1.519 &    0.210 &    0.527 \\
B &    5.165 &    0.103 &    0.355 &    0.451 &    1.247 &    0.228 &    0.285 \\
C &    4.355 &    0.084 &    0.330 &    0.183 &    0.926 &    0.458 &    0.356 \\
D &    5.021 &    0.120 &    0.406 &    0.100 &    0.826 &    1.202 &    0.491 \\
E &    4.669 &    0.064 &    0.976 &    0.207 &    1.788 &    0.309 &    0.546 \\
F &    8.603 &    0.052 &    0.526 &    0.133 &    1.124 &    0.390 &    0.468
\end{tabular}

}

The results of approach $\mathcal{A}$ can be found in Table~\ref{table:results_A}.
We observe, that the constant velocity values are slightly above the average values given in Table~\ref{table:pipe_data}.
This supports our claim from Section~\ref{sec:linearization} that overestimating the velocity is more favorable to reduce the pressure difference error.
The actual average error values are quite small in terms of absolute values.
However, when comparing to the actual average friction values, the values are quite high with an error to friction ratio above 20\%.
For pipe $D$ the error is even higher than the friction-based difference.
When looking at the maximum values, the relative errors are mainly at a level of about 50\%, which is also rather high.

In contrast to the results expected in Section~\ref{sec:real_data_analysis}, the values for pipes $A$ to $C$ are not significantly better than the ones of pipes $D$ to $F$.
The only clear difference between the pipes in the results is the bad average error of pipe $D$ in relation to the friction values.

{
\centering
\small
\captionsetup{justification=centering}
\captionof{table}{Results of approach $\mathcal{B}$ including the average pressure difference error, the corresponding maximum error, the average and maximum real friction-based pressure differences and the ratio of error to friction values for the average and maximum case.}
\label{table:results_B}
\begin{tabular}{lrrrrrr}
\vspace{0.05cm} Pipe & avg err & max err & avg $|fL|$ & max $|fL|$ & avg err/ & max err/\\
 & in bar & in bar & in bar & in bar & avg $|fL|$ & max $|fL|$\\
\midrule
A &    0.104 &    0.938 &    0.621 &    1.519 &    0.167 &    0.618 \\
B &    0.100 &    0.675 &    0.451 &    1.247 &    0.222 &    0.542 \\
C &    0.052 &    0.399 &    0.183 &    0.926 &    0.284 &    0.431 \\
D &    0.047 &    0.781 &    0.100 &    0.826 &    0.471 &    0.945 \\
E &    0.048 &    0.878 &    0.207 &    1.788 &    0.234 &    0.491 \\
F &    0.032 &    0.613 &    0.133 &    1.124 &    0.238 &    0.545
\end{tabular}

}

For approach $\mathcal{B}$ the results are shown in Table~\ref{table:results_B}.
For all pipes, the average values are smaller for approach $\mathcal{B}$.
The best improvement is made on pipe $D$, where the error could be more than halved.
However, in relation to the average friction, pipe $D$ has still by far the highest ratio.

For the maximum values, $\mathcal{B}$ performs even worse than $\mathcal{A}$ on all pipes except for pipe $E$.
On pipe $B$ and $D$ the error value nearly doubled, which leads to a maximum error nearly as big as the actual friction on pipe $D$.

Regarding the expected performance due to the investigated velocity changes in Section~\ref{sec:real_data_analysis}, we observe that the values for pipe $D$ are better than expected, but are clearly still the worst among all pipes, which is consistent with Figure~\ref{fig:velochanges}.
Despite the different values for hour 48 in the graphic, the other five pipes have quite similar results.
One reason could be the different evaluation periods: two years for Figure~\ref{fig:velochanges} and only the second year for the results of Table~\ref{table:results_B}.

\section{Conclusion}
We linearized the isothermal Euler Equations by fixing the velocity in the friction term to a constant value for each pipe.
The results of approach $\mathcal{B}$ for determining the constant value, where we fixed the velocity to the historic values of two days before, indicate that the average errors made are not too large compared to the overall friction-induced pressure drop on the pipelines.
However, the maximal error values turned out to be quite high, even in relation to the maximal friction values.
Hence, we can conclude that the presented fixed velocity approaches can only be considered as rough approximations.

For future research, the results on our six pipes have to be verified on the complete network.
Especially the bad values of pipe $D$ have to be analyzed to find potential structural problems, maybe due to the change in flow direction.
Furthermore, more sophisticated approaches to determine the fixed velocity might be possible.
One option is to examine if compressor configurations at network intersection points have an impact on the velocity of the adjacent pipelines.

\section*{Acknowledgments}
The work for this article has been conducted within the Research Campus MODAL funded by the German Federal Ministry of Education and Research (BMBF) (fund number 05M14ZAM).

%
%

\bibliographystyle{abbrv}
\bibliography{shortpaper}

\end{document}